\documentclass{article}
\usepackage{amsthm,amssymb,amsmath}
\usepackage{epsfig}
\usepackage{color}
\usepackage{tikz}

\newtheoremstyle{thm1}
{15pt}
{15pt}
{}
{}
{\bf}
{}
{\newline }
{}
\theoremstyle{thm1}

\newtheorem{defn}{thm1}[section]
\newtheorem{thm}[defn]{Theorem}
\newtheorem{lem}[defn]{Lemma}

\newtheorem{exmp}[defn]{Example}

\newtheorem{obs}[defn]{Observation}

\newtheorem{prob}[defn]{Problem}
\newtheorem{def1}[defn]{Definition}

\newenvironment{prf}[1][ ]{\textbf{Proof.} #1\\}{\hfill\qedsymbol\\[15pt]}

\begin{document}
\begin{center}
{\Large Graphs with degree complete labeling}\\
\vspace{3mm}
{ Sebastian Milz}\\
{Lehrstuhl II f\"{u}r Mathematik\\ RWTH Aachen University\\ 52062 Aachen, Germany}\\
{\tt milz@math2.rwth-aachen.de}\\
\end{center}
\begin{abstract}
  In 2006 Qian \cite{Qian06} introduced the concept of degree complete graphs for labeled graphs.
  He also gave a characterization of these graphs in terms of two forbidden subgraphs.
  Furthermore, he mentioned that the property of being degree complete depends on the labeling of the graph. 
  Related to this he stated the problem to find a characterization of those (unlabeled) graphs for which every labeled version is not degree complete.\\
  We say that a (unlabeled) graph has a degree complete labeling, if there is a labeled version of the graph that is degree complete.
  In this paper we give three characterizations of graphs with degree complete labeling.
  These characterizations give us polynomial-time procedures to recognize these graphs and find a degree complete labeling, if it exists.\\
  
\end{abstract}
 \section{Terminology and Introduction}
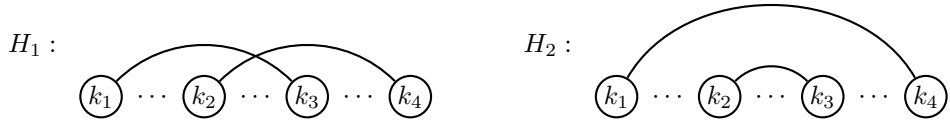
\begin{figure}[ht]\label{fig_H1H2}
  \begin{center}
   \begin{tikzpicture}[scale=0.65, vertex/.style={circle,inner sep=1pt,draw,thick}, myarrow/.style={thick}]
         \node at (-40pt,30pt) {$H_1:$};
         \node at (31pt,0pt) {$\cdots$};
         \node at (91pt,0pt) {$\cdots$};
         \node at (151pt,0pt) {$\cdots$};

         \node (1) at (0pt,0pt) [vertex] {$k_1$};
         \node (2) at (60pt,0pt) [vertex] {$k_2$};
         \node (3) at (120pt,0pt) [vertex] {$k_3$};
         \node (4) at (180pt,0pt) [vertex] {$k_4$};

         \draw[myarrow] (1) to[bend left=45]  (3);
         \draw[myarrow] (2) to[bend left=45] (4);

         \node at (260pt,30pt) {$H_2:$};
         \node at (331pt,0pt) {$\cdots$};
         \node at (391pt,0pt) {$\cdots$};
         \node at (451pt,0pt) {$\cdots$};
        
         \node (1a) at (300pt,0pt) [vertex] {$k_1$};
         \node (2a) at (360pt,0pt) [vertex] {$k_2$};
         \node (3a) at (420pt,0pt) [vertex] {$k_3$};
         \node (4a) at (480pt,0pt) [vertex] {$k_4$};

         \draw[myarrow] (1a) to[bend left=60]  (4a);
         \draw[myarrow] (2a) to[bend left=45] (3a);        
 \end{tikzpicture}
\end{center}
\caption{Forbidden configurations $H_1$ and $H_2$ from Theorem \ref{thm_Qian}.}
\end{figure}

 Let $G$ be a finite undirected graph with multiple edges but without loops.
 We denote $V(G)$ as the vertex set and $E(G)$ as the edge set of $G$.
 By $n=n(G)=|V(G)|$ and $m=m(G)=|E(G)|$ we refer to the {\it order} and the {\it size} of $G$, respectively.
 For an edge $e\in E(G)$ we use the notations $uv$ and $\{u,v\}$, if $e$ connects the vertices $u$ and $v$ in $G$. 
 Moreover, we says that $u$ and $v$ are {\it adjacent} and $u$ is a {\it neighbor} of $v$.
 For a vertex $v\in V(G)$ denote $N_G(v)$ the set of neighbors and $d_G(v)=|N_G(v)|$ the {\it degree} of $v$.\\ 
 A {\it subgraph} $H$ of $G$ is a graph satisfying $V(H)\subseteq V(G)$ and $E(H)\subseteq E(G)$.
 For any subset of edges $F\subseteq E(G)$ we define $G-F$ as the subgraph with vertex set $V(G)$ and edge set $E(G)\setminus F$.
 Considering a vertex set $X\subseteq V(G)$ a subgraph $H$ of $G$ is {\it induced} by $X$ and write $H=G[X]$, if $V(H)=X$ and $E(H)=\{uv \in E(G)\,|\, u,v\in X\}$.
 Furthermore, we use the notation $G-X$ for the induced subgraph of $G$ with vertex set $V\setminus X$.
 We simply write $G-x$ instead of $G-\{x\}$.\\
 Let $\{v_1,\ldots,v_{p+1}\}\subseteq V(G)$ and $\{e_1\ldots,e_{p}\}\subseteq E(G)$. 
 The sequence $v_1e_1v_2e_2\ldots e_{p}v_{p+1}$ is a {\it path}, if the vertices $v_1,\ldots,v_{p+1}$ are pairwise distinct and $e_i=v_iv_{i+1}$ for $1\le i\le p$.
 The same sequence is a {\it cycle}, if $v_1,\ldots,v_{p}$ are pairwise distinct, $v_1=v_{p+1}$, and $e_i=v_iv_{i+1}$ for $1\le i\le p$. 
 In both cases $p$ is the length of the path or cycle, respectively.\\ 
 We say $G$ is {\it connected}, if there is a path from $u$ to $v$ for every pair of distinct vertices $u,v\in V(G)$.
 By a {\it component} we refer to a maximal connected induced subgraph of $G$.
 A tree is a connected graph that does not contain any cycle of positive length.\\
 A {\it labeling} of a graph $G$ is a bijective function $f:V(G)\rightarrow \{1,\ldots,n\}$.
 The graph $G$ together with a labeling $f$ is a {\it labeled graph} denoted by $G_f$. 
 For simplicity we usually omit $f$ and identify the vertex set of a labeled graph with $\{1,\ldots,n\}$.\\
 An {\it orientation} of $G$ is a directed graph $D$ with vertex set $V(D)=V(G)$ in which every edge $uv\in E$ is assigned with a direction. 
 Thus, a directed edge is an ordered pair of vertices that we denote as an {\it arc}.
 We write $A(D)$ for the set of arcs in $D$ and $(u,v)$ for an element of $A(D)$, if the direction of the arc in $D$ is from $u$ to $v$.
 In this case we also denote $v$ as an {\it out neighbor} of $u$.
 By the {\it out degree} of $v$ we refer to the total number of out neighbors of a given vertex $v$ and write $d^+_D(v)$.
 Obviously, we have $0\le d^+_D(v)\le d_G(v)$ for every $v\in V(G)$.\\ 
 Now, let $G$ be a labeled graph with vertex set $\{1,2,\ldots,n\}$ and $D$ an orientation of $G$.
 A vector of nonnegative integers $s=(s_1,\ldots,s_n)$ is the {\it out degree vector} of $D$ or just a {\it degree vector} of $G$, if $d^+_D(i)=s_i$ for $1\le i\le n$. 
 In this case we also say that the out degree vector $s$ {\it is realized} by $D$.\\
 The following partial order ``$\preccurlyeq$'' on the set of nonnegative integer vectors is closely related to the well-known dominance order. 
 Let $s$ and $t$ be two nonnegative integer vectors of dimension $n$, then $$s\preccurlyeq t:\Leftrightarrow \sum_{i=1}^k{s_i}\le \sum_{i=1}^k{t_i}\quad \textrm{for all } 1\le k\le n$$ and equality holds for $k=n$.\\[3mm]
 Consider the two orientations $D^r$ and  $D^l$ of $G$ with arc sets $A_G^r=\{(i,j)\,|\, ij \in E(G), i<j\}$ and $A_G^l=\{(i,j)\,|\, ij \in E(G), i>j\}$.
 We define $s_G^r$ and $s_G^l$ as the out degree vectors $D^r$ and $D^l$, respectively.
 To visualize these orientations suppose all vertices of $G$ are written on a horizontal line with increasing labels from left to right.
 Now, $A_G^r$ is the arc set of the orientation of $G$, where all edges are oriented from left to right. 
 In the same arrangement of vertices all arcs in $A_G^l$ point from right to left.
 Hence we have 
 $$\sum_{i=1}^k{\left(s_G^l\right)_i}=m\left(G[\{1,\ldots,k\}]\right)\quad \text{and}\quad \sum_{i=1}^k{\left(s_G^r\right)_i}=m\left(G[\{1,\ldots,k\}]\right)+m_k,$$
 for $1\le k\le n$, where $m_k$ is the number of edges between the vertex sets $\{1,\ldots,k\}$ and $\{k+1,\ldots,n\}$ in $G$.
 Thus, on the one hand, for every out degree vector $s$ of $G$ holds
 \begin{align}
  s_G^l \preccurlyeq s \preccurlyeq s_G^r\quad \textrm{and}\quad 0\le s_i\le d_G(i)\;\textrm{ for }\; i=1,\ldots,n. \label{PartialOrderCond}
 \end{align}
 On the other hand, there may be nonnegative integer vectors satisfying (\ref{PartialOrderCond}) which are not realized by an orientation of $G$.\\
 In \cite{Qian06} Qian introduced the concept of degree complete graphs.
 A labeled graph $G$ is {\it degree complete}, if every nonnegative integer vector $s$ satisfying (\ref{PartialOrderCond}) is a degree vector of $G$. 
 Qian also proved the following characterization of degree complete graphs.
 \begin{thm}[Qian \cite{Qian06} (2006)]\label{thm_Qian}
 A labeled graph $G$ is degree complete if and only if $G$ does not contain one of the two subgraphs $H_1$ and $H_2$, where
 $$V(H_1)=V(H_2)=\{k_1,k_2,k_3,k_4\},\quad k_1<k_2<k_3<k_4$$
 and
 $$E(H_1)=\{k_1k_3,k_2k_4\},\quad E(H_2)=\{k_1k_4,k_2k_3\}.$$
 \end{thm}
 As in \cite{Qian06}, we will denote the subgraphs $H_1$ and $H_2$ also as {\it forbidden configurations}.\\
 Reminding the embedding of $G$ along a horizontal line with respect to the vertex indices these forbidden configurations can be visualized (see Figure \ref{fig_H1H2}) in the following way.
If we draw all edges on one site (above or below) of the line, the subgraph $H_1$ is a pair crossing independent edges. 
 Similarly, $H_2$ refers to a pair of overleaping independent edges in $G$.\\
An important property concerning the concept of degree complete graphs is illustrated by the following example.
 \begin{exmp}[Qian \cite{Qian06} (2006)]\label{exmp1}
Consider the labeled graphs $G_1$ and $G_2$ from Figure \ref{fig_exmp1}.
 \begin{figure}[ht]
 \begin{center}
 \begin{tabular}{cp{1cm}c}
  \begin{tikzpicture}[scale=0.7, vertex/.style={circle,inner sep=2pt,draw,thick}, myarrow/.style={thick}]
         \node at (-40pt,40pt) {$G_1$:};
         \node (1) at (0pt,0pt) [vertex] {$1$};
         \node (2) at (60pt,0pt) [vertex] {$2$};
         \node (3) at (120pt,0pt) [vertex] {$3$};
         \node (4) at (180pt,0pt) [vertex] {$4$};
         
         \draw[myarrow] (1) to (2);
         \draw[myarrow] (2) to (3);
         \draw[myarrow] (3) to (4);        
 \end{tikzpicture} & &
  \begin{tikzpicture}[scale=0.7, vertex/.style={circle,inner sep=2pt,draw,thick}, myarrow/.style={thick}]
         \node at (-40pt,40pt) {$G_2$:};
         \node (1) at (0pt,0pt) [vertex] {$1$};
         \node (2) at (60pt,0pt) [vertex] {$2$};
         \node (3) at (120pt,0pt) [vertex] {$3$};
         \node (4) at (180pt,0pt) [vertex] {$4$};
         
         \draw[myarrow] (1) to[bend left=45] (3);
         \draw[myarrow] (2) to (3);
         \draw[myarrow] (2) to[bend left=45] (4);
\end{tikzpicture}
 \end{tabular}
 \end{center}
  \caption{Graphs $G_1$ and $G_2$ from Example \ref{exmp1}.}\label{fig_exmp1}
 \end{figure}
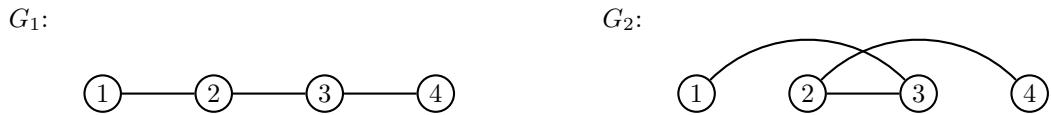
 Obviously, $G_1$ does not contain any of the subgraphs $H_1$ and $H_2$. 
 Thus, from Theorem \ref{thm_Qian} follows that $G_1$ is degree complete.
 On the other hand, in $G_2$ the edges $\{1,3\}$ and $\{2,4\}$ form a forbidden configuration $H_2$. 
 Therefore $G_2$ is not degree complete.
 \end{exmp}
 Since $G_1$ and $G_2$ are both labeled versions of a path of length $3$ we observe that the property of being degree complete depends on the vertex labeling of the graph.
 Qian also noticed this fact and stated the following problem.
 \begin{prob}[Qian \cite{Qian06} (2006)]\label{prob_Qian}
  Characterize the graphs which are not degree complete no matter how we label its vertices.
 \end{prob}
 To approach this topic we say that an unlabeled graph $G$ has a {\it degree complete labeling}, if there exists a labeling $f$ of its vertices such that the labeled graph $G_f$ is degree complete.
 Thus, Problem \ref{prob_Qian} asks for a characterization of the graphs which do not have any degree complete labeling.\\
 In addition to the above mentioned problem two more questions arise. 
 Firstly, can we find an efficient procedure to recognize graphs having a degree complete labeling?  
 Secondly, if we know that a given graph has such a labeling, how can we determine it?\\
 An unlabeled graph has a huge number of vertex labeling in general (even if we just count the labelings modulo the automorphism group of the graph). 
 Therefore, it is not useful to test different labelings by applying Theorem \ref{thm_Qian} and we need a new approach for this problem.\\
 The main theorem of the next section gives us three characterizations of unlabeled graphs which have a degree complete labeling. 
 The first equivalent formulation describes the structure of these graphs in terms of forbidden subgraphs.
 The other characterizations yield a polynomial procedure to recognize unlabeled graphs with degree complete labeling.
 Furthermore, these characterizations can be used as starting points for two similar algorithms which determine a desired labeling.
 \section{Degree complete labeling}
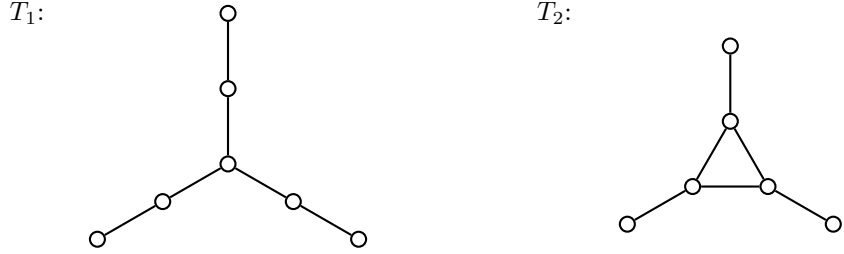
\begin{figure}[ht]
 \begin{center}
 \begin{tikzpicture}[scale=0.95, vertex/.style={circle,inner sep=2pt,draw,thick}, myarrow/.style={thick}]
       
         \node at (-80pt,60pt) {$T_1$:};
         \node (1) at (0pt,0pt) [vertex] {};
         \node (2) at (0pt,30pt) [vertex] {};
         \node (3) at (26pt,-15pt) [vertex] {};
         \node (4) at (-26pt,-15pt) [vertex] {};
         \node (5) at (0pt,60pt) [vertex] {};
         \node (6) at (52pt,-30pt) [vertex] {};
         \node (7) at (-52pt,-30pt) [vertex] {};
         
         \draw[myarrow] (1) to (2);
         \draw[myarrow] (1) to (3);
         \draw[myarrow] (1) to (4);
         \draw[myarrow] (2) to (5);
         \draw[myarrow] (3) to (6);
         \draw[myarrow] (4) to (7);

        \node at (130pt,60pt) {$T_2$:};
        \node (v1) at (200pt,17pt) [vertex] {};
        \node (v2) at (215pt,-9pt) [vertex] {};
        \node (v3) at (185pt,-9pt) [vertex] {};
        \node (v4) at (200pt,47pt) [vertex] {};
        \node (v5) at (241pt,-24pt) [vertex] {};
        \node (v6) at (159pt,-24pt) [vertex] {};
        
        \draw[myarrow] (v1) to (v2);
        \draw[myarrow] (v1) to (v3);
        \draw[myarrow] (v2) to (v3);
        \draw[myarrow] (v1) to (v4);
        \draw[myarrow] (v2) to (v5);
        \draw[myarrow] (v3) to (v6);
        
 \end{tikzpicture}
\end{center}
   \caption{Graphs $T_1$ and $T_2$}\label{fig_T1T2}
 \end{figure}

 We start with a simple but important observation which is a direct consequence of Theorem \ref{thm_Qian}.
 Moreover it motivates a characterization of graphs with degree complete labeling which is based on forbidden subgraphs.
 \begin{obs}\label{obs_subgraph}
  Let $G$ be a graph and $H$ a subgraph of $G$. If $G$ has a degree complete labeling, then $H$ has a degree complete labeling.
 \end{obs}
 Next, we want to obtain some necessary conditions for unlabeled graphs with degree complete labeling.
 Thus we consider those graphs that contain at least one of the forbidden configurations $H_1$ and $H_2$ in every labeling.
 In particular, we are interested in a set of pairwise not including graphs with this property.
 Denote $C_k$ the graph consisting of a cycle of length $k\ge 3$. 
 Furthermore, we define the graphs $T_1$ and $T_2$ as in Figure \ref{fig_T1T2}. 
 \begin{lem}\label{lem1}
 Let $G\in \{T_1,T_2\}\, \cup \, \{C_k\,|\, k\ge 4\}$. For every vertex labeling $f$ of $G$ there exists a forbidden configuration $H_1$ or $H_2$ in $G_f$.
 \end{lem}
 \begin{prf}
  We consider an arbitrary vertex labeling of $G$.
  Therefore, without loss of generality we identify each element in $V(G)$ with exactly one of the integers from 1 to $|V(G)|$.
  There are three cases.\\[2mm]
  {\it Case 1: $G=T_1$.}\\ 
  Denote $i_1$ the unique vertex of degree 3. 
  There are three distinct vertices in $V(G)$ which are adjacent to $i_1$. 
  Furthermore, we can find two of these vertices $i_2$ and $i_3$ such that exactly one the following two conditions holds.
  Either we have $i_1<i_2<i_3$ or $i_3<i_2<i_1$. 
  In both cases there is a vertex $i_4\in V(G)$ that has $i_2$ as its unique neighbor.  
  If $i_1<i_4<i_3$ (respectively $i_3<i_4<i_1$), then $\{i_1i_3,i_2i_4\}$ forms a forbidden configuration $H_2$. 
  Otherwise, $\{i_1i_3,i_2i_4\}$ is the edge set of a copy of $H_1$.\\[2mm]
  {\it Case 2: $G=T_2$.}\\
  Denote $i_1,i_2,i_3$ the three vertices of degree 3 in $T_2$ such that $i_1<i_2<i_3$. 
  There is a vertex $i_4\in V(G)$ that has $i_2$ as unique neighbor. 
  Now, we have to distinct two cases.
  If either $i_1<i_4<i_2$ or $i_2<i_4<i_3$ holds, then $\{i_1i_3,i_2i_4\}$ is the edge set of a forbidden configuration $H_2$ in $G$. 
  Otherwise we have $i_4<i_1$ or $i_4>i_3$. 
  In this case $G$ contains the edges $i_1i_3$ and $i_2i_4$, that is a copy of $H_1$.\\[2mm]
  {\it Case 3: $G=C_k$, $k\ge 4$.}\\
  We observe that vertex $1$ has a neighbor $i_1$ with $2<i_1\le k$. 
  Moreover, vertex $2$ is adjacent to a vertex $i_2$ satisfying $2<i_2\le k$ and $i_1\neq i_2$. 
  If $i_1<i_2$, then the edges $\{1,i_1\}$ and $\{2,i_2\}$ form a forbidden configuration $H_1$. 
  In the case $i_1>i_2$ the same edges give us $H_2$ as a subgraph of $G$.
 \end{prf}
 Combining Theorem \ref{thm_Qian}, Observation \ref{obs_subgraph}, and Lemma \ref{lem1} we deduce that every graph containing a subgraph $T_1,\,T_2$ or a cycle of length $k\ge4$ does not have a degree complete labeling.
 This shows that graphs with degree complete labeling have a structure that is similar to trees since they may only have cycles of length 3.\\[3mm]
 Now suppose for a moment that $G$ is a tree. 
 Obviously, $G$ cannot contain a subgraph isomorphic to $T_2$ or a cycle but it can have copy of $T_1$ as a subgraph.
 A tree without a copy of $T_1$ has a path $P$ such that all vertices which are not included in this path have degree 1 and are adjacent to a vertex of $P$.
 These trees are also known as caterpillars.
 A {\it Caterpillar} is defined (see \cite{Harary73}) as a tree $G$ such that $G-X_1$ is a path, where $X_1$ denotes the set of vertices of degree 1 in $G$.
 It is not difficult to see that caterpillars can be characterized as such trees without a subgraph isomorphic to $T_1$.
  
\begin{figure}
  \begin{center}
\begin{tikzpicture}[scale=0.8, vertex/.style={circle,inner sep=2pt,draw,thick}, myarrow/.style={thick}]

         \node (11) at (-161pt,0pt) [vertex] {};
         \node (12) at (-182pt,21pt) [vertex] {};
         \node (13) at (-182pt,-21pt) [vertex] {};
         \node (14) at (-131pt,0pt) [vertex] {};
         \node (15) at (-131pt,30pt) [vertex] {};
         \node (16) at (-101pt,0pt) [vertex] {};
         \node (17) at (-101pt,30pt) [vertex] {};
         \node (18) at (-80pt,21pt) [vertex] {};
         \node (19) at (-80pt,-21pt) [vertex] {};
         \node (20) at (-101pt,-30pt) [vertex] {};

         \draw[myarrow] (11) to (12);
         \draw[myarrow] (11) to (13);
         \draw[myarrow] (11) to (14);
         \draw[myarrow] (14) to (15);
         \draw[myarrow] (14) to (16);
         \draw[myarrow] (16) to (17);
         \draw[myarrow] (16) to (18);
         \draw[myarrow] (16) to (19);
         \draw[myarrow] (16) to (20);

         \node (1) at (0pt,0pt) [vertex,label=below:$1$] {};
         \node (2) at (30pt,0pt) [vertex,label=below:$2$] {};
         \node (3) at (60pt,0pt) [vertex,label=below:$3$] {};
         \node (4) at (90pt,0pt) [vertex,label=below:$4$] {};
         \node (5) at (120pt,0pt) [vertex,label=below:$5$] {};
         \node (6) at (150pt,0pt) [vertex,label=below:$6$] {};
         \node (7) at (180pt,0pt) [vertex,label=below:$7$] {};
         \node (8) at (210pt,0pt) [vertex,label=below:$8$] {};
         \node (9) at (240pt,0pt) [vertex,label=below:$9$] {};
         \node (10) at (270pt,0pt) [vertex,label=below:$10$] {};

         \draw[myarrow] (1) to (2);
         \draw[myarrow] (1) to[bend left=30] (3);
         \draw[myarrow] (1) to[bend left=40] (4);
         \draw[myarrow] (4) to (5);
         \draw[myarrow] (4) to[bend left=30] (6);
         \draw[myarrow] (6) to (7);
         \draw[myarrow] (6) to[bend left=30] (8);
         \draw[myarrow] (6) to[bend left=40] (9);
         \draw[myarrow] (6) to[bend left=50] (10);
   \end{tikzpicture}
  \end{center}
  \caption{An unlabeled caterpillar and labeled version of the same caterpillar with degree complete labeling.}\label{fig_caterpillar}
 \end{figure}
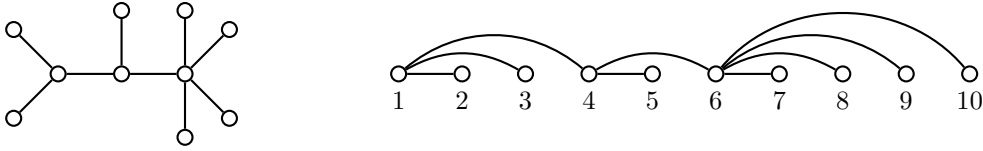

 Let $P=v_1e_1v_2e_2\ldots e_{p}v_{p+1}$, where $p$ is a nonnegative integer.
 Obviously, the vertex labeling $f$ defined by $f(v_i)=i$ is a degree complete labeling of $P$.
 Furthermore, we can extent $f$ to a degree complete labeling of the whole caterpillar $G$ by a repetition of the following step.
 Consider an unlabeled vertex $x\in X_1$ and denote $u\in V(P)$ the unique neighbor of $x$ in $G$.
 We add 1 to the label of every labeled vertex $w$ with $f(w)>f(u)$ and set $f(x)=f(u)+1$.
 Going on with this procedure we terminate with a labeling of all vertices that fulfills the following condition.
 For every edge $uv\in E(P)$ and every vertex $w$ satisfying $f(u)<f(w)<f(v)$ holds $w$ is adjacent to $u$.
 Hence $G_f$ does not contain a forbidden configuration $H_1$ or $H_2$.
 By Theorem \ref{thm_Qian} the labeling $f$ yields a degree complete labeling of $G$.\\
 On the left hand side of Figure \ref{fig_caterpillar} there is a caterpillar.
 The graph to the right refers to a labeled version of the same caterpillar. 
 The sequence of vertices of the labeled graph is with respect to the labeling and shows that it is a degree complete labeling.\\[3mm]
 We try to adapt the characterization known for caterpillars.
 Thus, we return to an arbitrary graph $G$ not containing $T_1,\,T_2$ as a subgraph or a cycle of length $k\ge4$. 
 Additional to the vertices of degree 1 in $X_1$ we have to delete a vertex or at least an edge of every triangle to obtain a path.
 Therefore we define the following sets.
 \begin{def1}\label{X1X2}
  Let $G$ be a graph with $V(G)=\{v_1,\ldots,v_n\}$ and $E(G)=\{e_1,\ldots,e_m\}$. We define the vertex set
 \begin{align*}
  X_1(G) &:=\{v\in V(G)\,|\, d_G(v)=1\}.
 \end{align*}
 Furthermore, we construct $X_2(G)$ by a procedure. Initialize $X_2(G)=\emptyset$. 
 For every $i$ from $1$ to $n$ we add $v_i$ to $X_2(G)$, if all following conditions are fulfilled
 \begin{itemize}
  \item $v_i$ has exactly two neighbors in $G$ denoted by $u$ and $v$.
  \item $v_i$ is the unique common neighbor of $u$ and $v$.
  \item $u$ and $v$ are adjacent.
  \item $u$ and $v$ are not in $X_2(G)$.
 \end{itemize}
 Finally, we define the edge set $F(G)$. 
 Start with $F(G)= \emptyset$.
 For every $j$ from $1$ to $m$ add the edge $e_j$ to $F(G)$, if all following conditions are fulfilled
 \begin{itemize}
  \item $e_j$ joints the vertices $u$ and $w$.
  \item $u$ and $w$ have a unique common neighbor $v$.
  \item $uv$ and $vw$ are not in $F(G)$.
 \end{itemize}
 \end{def1}
 By these definions we see that $w\in X_2(G)$ has the following properties.
 The vertex $w$ has degree 2 and is part of a triangle.
 Moreover, there does not exist a further triangle in $G$ containing both neighbors of $w$.
 We also notice that for every triangle in $G$ at most one of its vertices is a part of $X_2(G)$. 
 In general there is more than one possible choice for a maximal set satisfying all conditions of $X_2(G)$.\\
 Now, consider a triangle in $G$ where at least one vertex $v$ has degree 2. 
 By its definition $F(G)$ contains the edge of the triangle opposing $v$, if this edge is not contained in a further triangle.
 Similar to $X_2(G)$ there does not exist a triangle which has two edges in $F(G)$.\\
 The above mentioned constructions show the following. 
 The procedures which yield $X_2(G)$ and $F(G)$ are basically greedy algorithms that additionally use some graph search elements.
 Thus, $X_2(G)$ and $F(G)$ can be determined in polynomial time with respect to the order of $G$.\\[3mm]
 We are now able to formulate and prove the main theorem of this section.
 \begin{thm}\label{thm_char_list}
 Let $G$ be a graph. The following statements are equivalent:
 \begin{enumerate}
  \item[(i)] $G$ has a degree complete labeling.
  \item[(ii)] $G$ does not contain a subgraph isomorphic to $T_1$, $T_2$ or $C_k$ ($k\ge 4$).
  \item[(iii)] $G-X_1(G)- X_2(G)$ is a disjoint union of paths.
  \item[(iv)] $G-X_1(G)-F(G)$ is a disjoint union of paths.
 \end{enumerate}
\end{thm}
 \begin{prf}
  From (i) to (ii): 
  Follows from Theorem \ref{thm_Qian}, Observation \ref{obs_subgraph}, and Lemma \ref{lem1}.\\[3mm]
  From (ii) to (iii) and (ii) to (iv): 
  Suppose $G$ does not contain a subgraph isomorphic to $T_1$, $T_2$ or $C_k$ ($k\ge 4$).
  If $G$ is not connected we consider each component separately. 
  Thus we assume that $G$ is connected.
  There does not exist an edge in $G$ which is part of two triangles, otherwise $G$ contains a cycle of length 4. 
  Furthermore, every triangle has a vertex $v$ of degree 2 since $T_2$ is not a subgraph in $G$.
  Notice that $G-v$ is connected. 
  Hence $G-X_2(G)$ is connected because $X_2(G)$ contains at most one vertex in a triangle.
  On the other hand $X_2(G)$ has also at least one vertex in every triangle. 
  Thus, $G-X_2(G)$ is an induced tree in $G$. 
  Taking into account that $T_1$ is not a subgraph in $G-X_2(G)$ we deduce that it is a caterpillar.
  From the definition of caterpillars follows that $G-X_1(G)- X_2(G)$ is a path.\\
  Considering $F(G)$ we observe by the same arguments as before that every triangle of $G$ contains exactly one edge of $F(G)$. 
  Therefore, $G-F(G)$ is a connected subgraph of $G$ without a cycle.
  Hence it is a tree not including $T_1$, that is a caterpillar.
  Again, $G-X_1(G)-F(G)$ is a path.\\[3mm]
  From (iii) to (i):
  Let $G$ be a graph such that $G-X_1(G)- X_2(G)$ is a disjoint union of paths.
  Obviously, $G$ has a degree complete labeling if and only if every  component of $G$ has a degree complete labeling.
  Furthermore, if a component does not consist of a single edge, then it corresponds to exactly one of the paths in $G-X_1(G)- X_2(G)$.
  Since a single edge has a (trivial) degree complete labeling we can assume that there is a path in $G-X_1(G)- X_2(G)$ for each component.
  In the following we give a labeling procedure for an arbitrary  component of $G$ and prove that the obtained labeled graph is degree complete.\\
  Suppose $\tilde{G}$ is a component of $G$ and $P$ the corresponding path in $G-X_1(G)- X_2(G)$. 
  We define $$X_1(\tilde{G})=X_1(G)\cap V(\tilde{G})\quad \text{and}\quad X_2(\tilde{G})=X_2(G)\cap V(\tilde{G}).$$
  Let $P$ be of length $p\ge0$ and denote $v_1,\ldots,v_{p+1}$ the vertices of $P$ such that $v_i$ and $v_{i+1}$ are adjacent for $1\le i\le p$.
  First we initialize the labeling function $f$ by $f(v_i)=i$.
  Next, we extent $f$ to the vertices in $X_2(\tilde{G})$ one by one as follows.
  From the definition of $X_2(\tilde{G})$ we observe that every unlabeled vertex $v\in X_2(\tilde{G})$ has exactly two neighbors, say $v_j$ and $v_{j+1}$, in $P$ which satisfy $f(v_j)<f(v_{j+1})$.
  Now, we add 1 to the label of every labeled vertex $w$ with $f(w)>f(v_j)$ and set $f(v)=f(v_j)+1$.  
  If all elements in $X_2(\tilde{G})$ are labeled, we finally consider an unlabeled vertex in $x\in X_1(\tilde{G})$.
  Notice that $x$ has a labeled vertex $u\in V(P)$ as its unique neighbor.
  Again, we add 1 to the label of every labeled vertex $w$ with $f(w)>f(u)$ and set $f(x)=f(u)+1$.\\ 
  It is not difficult to see that $f$ is a bijective function from $V(\tilde{G})$ to $\{1,\ldots,|V(\tilde{G})|\}$, that is a labeling of $\tilde{G}$.
  To verify that $f$ is also a degree complete labeling, let $v_iv_{i+1}\in E(P)$ and consider an arbitrary vertex $y$ such that $f(v_i)<f(y)<f(v_{i+1})$.
  Obviously, we have $y\in X_1(\tilde{G})\cup X_2(\tilde{G})$.
  If $y\in X_2(G)$, then $y$ is adjacent to both $v_i$ and $v_{i+1}$.
  Furthermore, by its construction there does not exist a further vertex in $X_2(G)$ with the property of $y$.
  Thus every other vertex $z$ with $f(v_i)<f(z)<f(v_{i+1})$ is in $X_1(\tilde{G})$.
  From the above mentioned labeling procedure we deduce that $z$ has $v_i$ as its unique neighbor and $f(z)<f(y)$.
  Therefore, the subgraph induced by all vertices with labels from $f(v_i)$ to $f(v_{i+1})$ does not contain a forbidden configuration from Theorem \ref{thm_Qian}.
  Moreover, this also holds for $\tilde{G}_f$ because $\tilde{G}$ does not have an edge $uv\in E(\tilde{G})$ such that $f(u)<f(v_i)<f(v)$ for $1\le i\le p+1$
  Hence $G_f$ hence is degree complete.\\[3mm]
  From (iv) to (i):
  Let $G$ be a graph such that $G-X_1(G)-F(G)$ is a disjoint union of paths.
  By similar arguments as mentioned before it is sufficient to show that each  component of $G$ has a degree complete labeling.
  Again, we assume that every component corresponds to exactly one of the paths in $G-X_1(G)-F(G)$.\\
  Let $\tilde{G}$ be a component of $G$ and $P$ the corresponding path in $G-X_1(G)-F(G)$ with vertices $v_1,\ldots,v_{p+1}$ such that $v_i$ and $v_{i+1}$ are consecutive in $P$ for $1\le i\le p$.
  We define $$X_1(\tilde{G})=X_1(G)\cap V(\tilde{G})\quad \text{and}\quad F(\tilde{G})=F(G)\cap E(\tilde{G}).$$
  Furthermore, we initialize $f$ by $f(v_i)=i$ for every $1\le i\le p+1$.
  Now, consider an unlabeled vertex $x\in X_1(\tilde{G})$ and denote $u$ the unique neighbor of $x$ in $P$.
  Again, we extent $f$ by adding 1 to the label of every labeled vertex $w$ with $f(w)>f(u)$ and set $f(v)=f(u)+1$.\\
  Obviously, $f$ is a labeling of $\tilde{G}$. 
  Moreover, for every edge $v_iv_{i+1}\in E(P)$ there are two possibilities.
  If there is an edge in $F(\tilde{G})$ joining $v_{i-1}$ and $v_{i+1}$, then $f(v_{i+1})=f(v_{i})+1$.
  Thus there does not exist a vertex with label between $f(v_i)$ and $f(v_{i+1})$ is this case.
  Otherwise there might be a vertex $z$ with $f(v_i)<f(z)<f(v_{i+1})$.
  As seen before, we have $z\in X_1(\tilde{G})$ and $z$ has $v_i$ as its unique neighbor.
  Therefore, $\tilde{G}_f$ is degree complete as it does not contain a forbidden configuration $H_1$ or $H_2$.
 \end{prf}
 We finish this section with an example on the labeling procedures from the proof of Theorem \ref{thm_char_list}.
 Let $G$ be the graph form Figure \ref{fig_exmp2}.
 \begin{figure}[ht]
 \begin{center}
 \begin{tikzpicture}[scale=0.6, vertex/.style={circle,inner sep=0.5pt,draw,thick}, vertex2/.style={circle,inner sep=2pt,draw,thick}, myarrow/.style={thick}]
   \node (1) at (-15pt,36pt) [vertex2] {$v_{5}$};
   \node (2) at (8pt,83pt) [vertex2] {$v_2$};
   \node (3) at (-50pt,81pt) [vertex2] {$v_1$};
   \node (4) at (70pt,1pt) [vertex] {$v_{10}$};
   \node (5) at (40pt,43pt) [vertex2] {$v_6$};
   \node (6) at (-25pt,-5pt) [vertex2] {$v_9$};
   \node (7) at (-71pt,38pt) [vertex2] {$v_4$};
   \node (8) at (154pt,3pt) [vertex] {$v_{11}$};
   \node (9) at (105pt,31pt) [vertex2] {$v_7$};
   \node (10) at (145pt,53pt) [vertex2] {$v_8$};
   \node (11) at (47pt,93pt) [vertex2] {$v_3$};

   \draw[myarrow] (1) to (3);
   \draw[myarrow] (1) to (5);
   \draw[myarrow] (1) to (6);
   \draw[myarrow] (1) to (7);
   \draw[myarrow] (2) to (5);
   \draw[myarrow] (3) to (7);
   \draw[myarrow] (4) to (5);
   \draw[myarrow] (4) to (9);
   \draw[myarrow] (5) to (9);
   \draw[myarrow] (5) to (11);
   \draw[myarrow] (8) to (9);
   \draw[myarrow] (9) to (10);
 \end{tikzpicture}
 \end{center}
 \caption{Graph $G$}\label{fig_exmp2}
 \end{figure}
 First, we determine the following sets 
 $$X_1(G)=\{v_2, v_3, v_8, v_9, v_{11}\},\quad X_2(G)=\{v_1,v_{10}\}\quad \text{and}\quad F(G)=\{v_4v_5,v_6v_7\}.$$
 We observe that $G-X_1(G)-X_2(G)$ is a path with vertex set $\{v_4,v_5,v_6,v_7\}$.
 Similarly, $G-X_1(G)-F(G)$ consists of the path $v_4v_1v_5v_6v_{10}v_7$.
 Thus, $G$ has a degree complete labeling.\\
 Next, as described in the proof from (iii) to (i), we initialize the labeling $f$ by 
 $$f(v_4)=1,\quad f(v_5)=2,\quad f(v_6)=3,\quad\text{and}\quad f(v_7)=4.$$
 Since $v_1\in X_2(G)$ is unlabeled and $N_G(v_1)=\{v_4,v_5\}$ we set 
 $$f(v_1)=2,\quad f(v_5)=2+1=3,\quad f(v_6)=3+1=4,\quad\text{and}\quad f(v_7)=4+1=5.$$
 Analogously, for $v_{10}\in X_2(G)$ the procedure yields
 $$f(v_{10})=5\quad\text{and}\quad f(v_7)=5+1=6.$$
 For the unlabeled vertex $v_2\in X_1(G)$ we continue with $$f(v_{2})=5,\quad f(v_{10})=5+1=6,\quad\text{and}\quad f(v_{7})=6+1=7.$$

 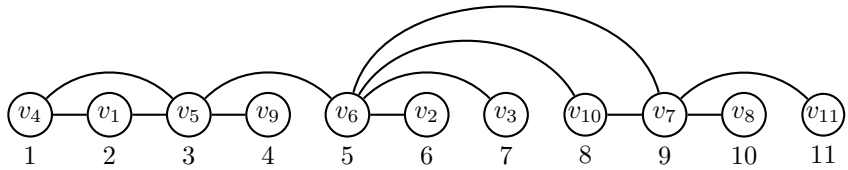
\begin{figure}[ht]
 \begin{center}
 \begin{tikzpicture}[scale=0.6, vertex/.style={circle,inner sep=0.5pt,draw,thick},
                     vertex2/.style={circle,inner sep=2pt,draw,thick},
                    myarrow/.style={thick},myarrow2/.style={thick,opacity=0.1}]
    \node (1) at (100pt,-70pt) [vertex2,label=below:$3$] {$v_5$};
   \node (2) at (300pt,-70pt) [vertex2,label=below:$7$] {$v_3$};
   \node (3) at (50pt,-70pt) [vertex2,label=below:$2$] {$v_1$};
   \node (4) at (350pt,-70pt) [vertex,label=below:$8$] {$v_{10}$};
   \node (5) at (200pt,-70pt) [vertex2,label=below:$5$] {$v_6$};
   \node (6) at (150pt,-70pt) [vertex2,label=below:$4$] {$v_9$};
   \node (7) at (0pt,-70pt) [vertex2,label=below:$1$] {$v_4$};
   \node (8) at (500pt,-70pt) [vertex,label=below:$11$] {$v_{11}$};
   \node (9) at (400pt,-70pt) [vertex2,label=below:$9$] {$v_7$};
   \node (10) at (450pt,-70pt) [vertex2,label=below:$10$] {$v_8$};
   \node (11) at (250pt,-70pt) [vertex2,label=below:$6$] {$v_2$};
   
   \draw[myarrow] (1) to (3);
   \draw[myarrow, bend left=45] (1) to (5);
   \draw[myarrow] (1) to (6);
   \draw[myarrow, bend right=45] (1) to (7);
   \draw[myarrow, bend right=45] (2) to (5);
   \draw[myarrow] (3) to (7);
   \draw[myarrow, bend right=60] (4) to (5);
   \draw[myarrow] (4) to (9);
   \draw[myarrow, bend left=75] (5) to (9);
   \draw[myarrow] (5) to (11);
   \draw[myarrow, bend right=45] (8) to (9);
   \draw[myarrow] (9) to (10);
 \end{tikzpicture}
 \end{center}
 \caption{Labeled Graph $G_f$}\label{fig_exmp2_labeled}
 \end{figure}
 Repeating this step until all vertices are labeled we arrive at the labeling $f$ from Figure \ref{fig_exmp2_labeled}.
 This also shows that $G_f$ is degree complete.
 Finally, it is not difficult to prove that the procedure from (iv) to (i) yields the same labeling.

\end{document}